\documentclass[11pt]{article}
\usepackage{latexsym}
\usepackage{comment}
\usepackage{amsmath,amsthm,amsfonts,amssymb,graphicx,epsfig,latexsym,color}
\usepackage{epsf,enumerate}

\title{A characterization of higher rank 
symmetric spaces via bounded cohomology}
  \author{Mladen Bestvina and Koji Fujiwara \thanks{The first author
  gratefully acknowledges the support by the National Science
  Foundation. The second author appreciates the hospitality of the
  mathematics department of the University of Utah.}}
  \date{\today}

\newtheorem{thm}{Theorem}[section]
\newtheorem{lemma}[thm]{Lemma}
\newtheorem{cor}[thm]{Corollary}
\newtheorem{prop}[thm]{Proposition}
\newtheorem*{main}{Main Theorem}{}
\newtheorem*{rrt}{Rank Rigidity Theorem}{}
{}

\theoremstyle{remark}

\newtheorem{definition}[thm]{Definition}
\newtheorem{remark}[thm]{Remark}

\newtheorem*{definition*}{Definition}
\newtheorem*{remark*}{Remark}

\def\R{{\mathbb R}}

\def\Z{{\mathbb Z}}

\def\H{{\mathbb H}}  
\def\G{{\Gamma}}

\def\WP{{\mathcal{WP}}}

\begin{document}

\maketitle
\begin{abstract} 
Let $M$ be complete nonpositively curved Riemannian manifold of finite
volume whose fundamental group $\Gamma$ does not contain a finite
index subgroup which is a product of infinite groups. We show that the
universal cover $\tilde M$ is a higher rank symmetric space iff 
$H^2_b(M;\R)\to H^2(M;\R)$ is injective (and otherwise the kernel is
infinite-dimensional). This is the converse of a theorem of
Burger-Monod. The proof uses the celebrated Rank Rigidity Theorem, as
well as a new construction of quasi-homomorphisms on groups that act
on $CAT(0)$ spaces and contain rank 1 elements.
\end{abstract}
\tableofcontents

\section{Introduction}\label{s:intro}

Let $\Gamma$ be a group. A function $\phi:\Gamma\to\R$ is called a {\it
  quasi-homomorphism} provided
$$\Delta(\phi):=\sup_{x,y\in\Gamma}|\phi(xy)-\phi(x)-\phi(y)|<\infty$$
The vector space
$$\widetilde{QH}(\Gamma)=\{\mbox{quasi-homomorphisms
}\phi:\Gamma\to\R\}/Hom(\Gamma,\R)+\{\mbox{bounded }\phi:\Gamma\to\R\}$$ 
can be identified with the kernel $Ker(H^2_b(\Gamma;\R)\to
  H^2(\Gamma;\R))$ from the bounded cohomology of $\Gamma$ to the
  ordinary cohomology.

A quasi-homomorphism $\phi$ is {\it homogeneous} if $\phi(g^n)=n\phi(g)$ for
all $g\in\Gamma$ and $n\in\Z$.
If we regard two quasi-homomorphisms $\phi_1,\phi_2$ as equivalent if
$\phi_1-\phi_2$ is bounded, then each equivalence class contains a unique
homogeneous representative. Note that homogeneous quasi-homomorphisms
are constant on conjugacy classes and that
$$\widetilde{QH}(\Gamma)=\{\mbox{homogeneous quasi-homomorphisms
}\phi:\Gamma\to\R\}/Hom(\Gamma,\R)$$

There is a growing list of groups $\Gamma$ for which
$\widetilde{QH}(\Gamma)$ has been ``computed''. If $\Gamma$ is
amenable then $\widetilde{QH}(\Gamma)=0$. The celebrated work of
Burger-Monod \cite{BuMo1}\cite{BuMo2} implies that
$\widetilde{QH}(\Gamma)=0$ for any irreducible lattice $\Gamma$ in a
semi-simple Lie group $G$ with finite center. (Burger and Monod
  state their result in the language of algebraic groups. In
  particular, their assumption that the group is simply-connected
  refers to the associated group over $\mathbb C$. In fact, it suffices that
  the complex group has finite fundamental group.) On the other hand,
non-elementary hyperbolic groups $\Gamma$, and indeed groups that act
on $\delta$-hyperbolic spaces with a weak proper discontinuity
condition (WPD) have infinite-dimensional $\widetilde{QH}(\Gamma)$
\cite{ef},\cite{f1},\cite{f2},\cite{bf}. That some condition on the
action is needed can be seen by looking at an irreducible lattice
$\Gamma$ in $\H^2\times\H^2$, which acts on $\H^2$ and has
$\widetilde{QH}(\Gamma)=0$ by the Burger-Monod theorem.

The constructions of quasi-homomorphisms in all these cases build on
the construction of Brooks \cite{b} for the case of free groups. In
this paper we further extend the technology to the setting of $CAT(0)$
spaces. More precisely, the main theorem is the following. For
definitions of rank 1 isometries and the WPD condition, see Sections
\ref{rk1} and \ref{wpd}.

\begin{main}
Let $X$ be a $CAT(0)$ space and let $\Gamma$ be a group acting by
isometries on $X$. Also suppose that $\Gamma$ is not virtually cyclic,
that the action satisfies $WPD$, and that $\Gamma$ contains an element
that acts as a rank 1 isometry. Then $\widetilde{QH}(\Gamma)$ is
infinite-dimensional.
\end{main}

In the Main Theorem we do not assume that $X$ is a proper, or even
complete, metric space. 

As an application, we prove the converse of the Burger-Monod
result. The proof is based on the celebrated Rank Rigidity Theorem
\cite[Theorem C, Theorem III.3.4]{bbook}, \cite[Theorem
  9.4.1]{Eb}.

\begin{rrt}
Let $M$ be a complete Riemannian manifold of nonpositive curvature and
finite volume. Consider the de Rham decomposition 
of the universal cover $\tilde M$ and
assume that it has no Euclidean factor. Then one of the following
holds.
\begin{enumerate}[(1)]
\item $\tilde M$ is a
symmetric space of noncompact type and rank $>1$.
\item Some deck transformation of $\tilde M$ is a
rank 1 isometry.
\item Some finite cover $M'$ of $M$ splits as a nontrivial Riemannian product
  $M'_1\times M'_2$.
\end{enumerate}
\end{rrt}

Our converse to the Burger-Monod theorem follows quickly.

\begin{thm}\label{converse}
Let $M$ be a complete Riemannian manifold of nonpositive curvature and
finite volume. Assume that $\Gamma=\pi_1(M)$ is finitely generated and
does not contain a
subgroup of finite index which is cyclic or a Cartesian product of two
infinite groups. Then the universal cover $\tilde M$ is a higher rank
symmetric space if and only if $\widetilde{QH}(\Gamma)=0$. Otherwise,
$\widetilde{QH}(\Gamma)$ is infinite-dimensional.
\end{thm}
\begin{proof}
Consider the de Rham decomposition of $\tilde M$. One
possibility is that $\tilde M$ is Euclidean space. Then $\Gamma$ is a
Bieberbach group and therefore contains $\Z^n$ as a subgroup of finite
index, contradicting our assumptions. If the de Rham decomposition
contains a Euclidean factor and a non-Euclidean factor, then some
finite cover $M'$ of $M$ is diffeomorphic to the product $N\times T$
for a nonpositively curved manifold $N$ of finite volume and a torus
$T$ \cite[Theorem 1.9]{BaEb} (this uses that the group is finitely
generated), again violating the assumptions. If there are no Euclidean
factors, then we apply the Rank Rigidity Theorem. Possibility (3) is
excluded by our assumption on $\Gamma$. 
If $\tilde M$
is a symmetric space of rank $>1$ the Burger-Monod theorem implies
$\widetilde{QH}(\Gamma)=0$, and if $\tilde M$ admits a rank 1 isometry
the Main Theorem implies that
$\widetilde{QH}(\Gamma)$ is infinite-dimensional.
\end{proof}

The condition $\widetilde{QH}(\Gamma)=0$ can be cast, by the work of
Bavard \cite{bavard}, in terms of stable commutator length. Recall
that for $\gamma\in [\Gamma,\Gamma]$ the {\it commutator length}
$c(\gamma)$ is the smallest $k\geq 0$ such that $\gamma$ can be
expressed as a product of $k$ commutators. The {\it stable commutator
length} is the limit $$c_\infty(\gamma)=\lim_{n\to\infty}\frac
{c(\gamma^n)}n$$
Bavard's theorem states that $\widetilde{QH}(\Gamma)=0$ if and only if
$c_\infty$ vanishes on $[\Gamma,\Gamma]$. 

There are examples of groups $\Gamma$ such that
$\widetilde{QH}(\Gamma)$ is nonzero and finite dimensional. All known
examples are constructed via central extensions as a variation of the
following argument (see \cite{MR2197372},\cite{zhuang}). Start with a
group $\Gamma_0$ with $\widetilde{QH}(\Gamma_0)=0$ and with
$H^2(\Gamma_0;\R)\neq 0$. Let $1\to\Z\to\Gamma\to\Gamma_0\to 1$ be a
central extension whose Euler class $e\in H^2(\Gamma_0;\R)$ is nonzero
and in the image of $H^2_b(\Gamma;\R)\to H^2(\Gamma;\R)$. Then
$H^2(\Gamma;\R)\cong H^2(\Gamma_0;\R)/<e>$ while
$H^2_b(\Gamma;\R)\cong H^2_b(\Gamma_0;\R)$, so
$\widetilde{QH}(\Gamma)$ is 1-dimensional.

\begin{thm}\label{dichotomy}
Suppose $M$ is a complete Riemannian manifold of nonpositive curvature
and finite volume. If $\Gamma=\pi_1(M)$ is finitely generated then
$\widetilde{QH}(\Gamma)$ is either 0 or infinite-dimensional.
\end{thm}

We will prove Theorem \ref{dichotomy} in Section \ref{reducible}. In
the special case when $\tilde M$ does not have Euclidean space as a de
Rham factor, and no finite cover of $M$ is a nontrivial Riemannian
product, the statement follows from Theorem \ref{converse}.

Another application is to the action of 
the mapping class group, $MCG(S)$, of a compact
orientable surface $S$ on the Teichm\"uller space, $T(S)$,
with Weil-Petersson metric.
We will show the following in Section \ref{wp}.
This result was shown in \cite{bf} using 
the action of $MCG(S)$ on the curve complex.

\begin{thm}\label{mapping class group}
Let $S$ be a compact orientable surface,
and $MCG(S)$ its mapping class group.
Suppose $\Gamma<MCG(S)$ is not virtually cyclic and 
contains at least one pseudo-Anosov element.
Then, $\widetilde{QH}(\Gamma)$
is infinite-dimensional.
\end{thm}

The present proof raises hopes that a similar theorem can be proved
for subgroups of $Out(F_n)$, since in that case there is no analog of
the curve complex (which is presently known to be hyperbolic). Of
course, there are still formidable obstacles, not the least of which
is finding a suitable metric on Outer space.

With an eye towards Outer space, we have tried to axiomatize the
properties of the space $X$ our method requires to show that
$\widetilde{QH}(\Gamma)$ is infinite dimensional. There are two
axioms: (DD) and (FT) (see the next section). Both hold in $CAT(0)$
spaces as well as in $\delta$-hyperbolic spaces.

We list the two axioms in Section \ref{s:axioms}. In Section
\ref{s:bcontr} we introduce the key concept of $B$-contracting
segments and discuss the basic properties. In Section \ref{s:qh} we
construct many quasi-homomorphisms in the situation when our axioms
hold. The key is a certain ``thick-thin'' decomposition of a geodesic
triangle (see the proof of Theorem \ref{Bcont.qh}). In Section
\ref{rk1} we introduce the concept of a rank 1 isometry and construct
a ``Schottky group'' (a free group of isometries in which every
nontrivial element is rank 1) starting with two independent rank 1
isometries. This is then used in Section \ref{wpd} to show that the
quasi-homomorphisms constructed in Section \ref{s:qh} span an
infinite-dimensional vector space. Finally, Section \ref{wp} contains
the verification of our axioms for the Weil-Petersson metric on
Teichm\" uller space.

\vskip 0.3cm
\noindent {\bf Acknowledgements.} We wish to thank Ken Bromberg for
pointing out that we should subdivide the ruled rectangles into two
triangles in the proofs of Section \ref{wp}. We also thank Werner
Ballmann, Emmanuel Breuillard, Ursula Hamenst\"adt and Nicolas Monod
for useful conversations and to Werner Ballmann for correcting a
mistake in an earlier version of the manuscript. We are grateful to
the referee for carefully reading the manuscript and suggesting we add
Section \ref{reducible}.

\section{Axioms}\label{s:axioms}
Let $(X,d)$ be a geodesic metric space. Thus any two points $a,b\in X$
are connected by at least one geodesic. To simplify notation we will
usually denote such a geodesic by $[a,b]$ even though it may not be
unique. We will also usually denote by $|a-b|$ the distance $d(a,b)$. 
For the rest of the paper we will make the following
assumptions on $X$.

For any geodesic segment $[a,b]$ and $x\in X$ denote by $\pi_{ab}(x)$
the (closed) set of points in $[a,b]$ that minimize the distance to
$x$. We think of $\pi_{ab}(x)$ as the (multivalued) projection of $x$
to $[a,b]$. The constant $C>0$ is fixed.

The first axiom states that the projection coarsely decreases
distances.

\begin{enumerate}
\item [{\bf (DD)}] For every $p\in\pi_{ab}(x)$ and every $p'\in
  \pi_{ab}(x')$ we have $$|p-p'|<|x-x'|+C.$$
\end{enumerate}

Note that, in particular, the diameter of $\pi_{ab}(x)$ is $<C$.

Our second axiom on $(X,d)$ is the strong version of the Fellow Traveller
Property.

\begin{enumerate}
\item [{\bf (FT)}] 
Suppose $|a-a'|\le D$, $|b-b'|\le D$. Then $[a',b']$ is contained in
the Hausdorff $C+D$-neighborhood of $[a,b]$.
\end{enumerate}

In particular, any two geodesics from $a$ to $b$ are in each others' Hausdorff
$C$-neighborhood.

\begin{prop}\label{DDFT}
Every $CAT(0)$ space (or more generally, a geodesic metric space with
a strictly convex distance function) as well as every
$\delta$-hyperbolic space satisfies axioms (DD) and (FT).\qed
\end{prop}

\begin{remark}
The strong version of the Fellow Traveller Property (FT) is used
only in the proof of Proposition \ref{Hausdorff neighborhoods}. For
the rest of the paper a weaker version (FT-) would suffice: $[a',b']$
is contained in the Hausdorff $FT(C,D)$-neighborhood for some function
$FT:[0,\infty)\times [0,\infty)\to [0,\infty)$. 
\end{remark}

\section{$B$-contracting segments}\label{s:bcontr}

From now on, we assume that $X$ is a CAT(0) space satisfying Axioms
(DD) and (FT).

\begin{definition}
A geodesic segment $[a,b]$ is said to be {\it $B$-contracting} for
$B>0$ if for every metric ball $K$ disjoint from $[a,b]$ the
projection $\pi_{ab}(K)$ has diameter $<B$.
\end{definition}

\begin{lemma}[Subsegments contracting]\label{subsegments contracting}
There is $B'>0$ (e.g. $B'=B+4C+3$ works) such that every subsegment of
a $B$-contracting segment is $B'$-contracting.
\end{lemma}

\begin{proof} Let $[a,b]$ be a $B$-contracting segment and
let $K$ be a ball centered at $O\in X$ and disjoint from $[u,v]\subset
[a,b]$. For concreteness, say $u$ is between $a$ and $v$. We will also
assume that $d(u,v)>C$ for otherwise $[u,v]$ is $C+1$-contracting.
There are two cases.

{\bf Case 1.} $\pi_{ab}(O)\cap [u,v]\neq\emptyset$. Then $K$ is also
disjoint from $[a,b]$ and so $\pi_{ab}(K)$ has diameter $<B$. Let
$x\in K$. If $\pi_{ab}(x)\cap [u,v]\neq\emptyset$ then
$\pi_{uv}(x)=\pi_{ab}(x)\cap [u,v]$ and if $\pi_{ab}(x)\cap
[u,v]=\emptyset$ then $\pi_{ab}(x)$ is contained either in $[a,u]$ or
in $[v,b]$ by Axiom (DD). Say $\pi_{ab}(x)\subset [a,u]$. Again by
(DD), it follows that $\pi_{uv}(x)$ is contained in the
$C$-neighborhood of $u$. (For suppose that $t\in \pi_{uv}(x)$ and
$|t-u|>C$. Let $s\in [a,u]$ be the last point with $|s-x|=|t-x|$. Then
both endpoints of $[s,t]$ are in $\pi_{st}(x)$ and $|s-t|>C$.) If
there are any such $x\in K$ then $u$ is in the smallest interval containing
$\pi_{ab}(K)$, so it follows that in either case the diameter of
$\pi_{uv}(K)$ is bounded above by $B+C$.

{\bf Case 2.} $\pi_{ab}(O)\cap [u,v]=\emptyset$. Again, for
concreteness, we may assume $\pi_{ab}(O)\subset [a,u]$. In this case
$K$ may intersect $[a,b]$, so we construct a different ball. 

\begin{figure}[htbp]
\begin{center}
\input{contracting-fig.pstex_t}
\caption{Lemma \ref{subsegments contracting}}
\end{center}
\end{figure}

Suppose
that there is $x\in
K$ such that $\pi_{uv}(x)$ is not contained in the
$B+4C+3$-neighborhood of $u$. It follows that $\pi_{ab}(x)$ is not
contained in the $B+3C+3$-neighborhood of $[a,u]$. By subdividing a
segment $[O,x]$ into intervals of size $<1$ and using Axiom (DD) we
find a point $y\in [O,x]$ such that $\pi_{ab}(y)$ is contained in the
$C+1$-neighborhood of $u$. Let $K'$ be the closed ball centered at $y$ with
radius $d(y,[a,b])-1$. Thus by assumption $\pi_{ab}(K')$ has diameter
less than $B$. 
Note that
$$|y-u|\geq |O-u|-|O-y|=|O-u|-|O-x|+|y-x|\geq |y-x|$$
(The last inequality follows from the fact that $u$ is outside $K$
while $x$ is inside.)

Now radius of $K'$ is $d(y,[a,b])-1\geq |y-u|-C-2$ and therefore
$K'$, even if it does not contain $x$, must contain a point at
distance $\leq C+2$ from $x$. The projection of this point to $[a,b]$
is then, by (DD), not contained in the
$B+3C+3-(2C+2)=B+C+1$-neighborhood of $[a,u]$. Since $\pi_{ab}(y)$ is
contained in the $C+1$-neighborhood of $u$ it follows that
$\pi_{ab}(K)$ has diameter $\geq B$, contradiction.
\end{proof}

For clarity, we will not try to keep up with constants such as
$B+4C+3$ in Lemma \ref{subsegments contracting}, 
but rather we introduce the following function
notation:
$\Phi=\Phi(B,C)$ stands for a positive function which is monotonically
increasing in each variable. 
Sometimes we decorate $\Phi$ with a subscript, and sometimes
there are more than 2
arguments. When we want to refer to a $\Phi$-function from a previous
lemma, we will indicate the number in the subscript,
e.g. $\Phi_{\ref{subsegments contracting}}(B,C)=B+4C+3$
and $\Phi_{\ref{thin}}(B,C)=3B+C+1$.

\begin{lemma}\label{5 balls}
Suppose $[x, y]$ and $[u, v]$ are geodesic segments such that
$u\in\pi_{uv}(x)$ and $v\in\pi_{uv}(y)$. Also suppose that the
distance between $[x, y]$ and $[u, v]$ is larger than $|u - v|$. Then
$[x, y]$ is contained in the union of 5 balls disjoint from $[u, v]$.
\end{lemma}

\begin{proof} 
Let $[p, q]$ be a shortest geodesic segment connecting a point $p\in
[x, y]$ with a point $q\in [u,v]$. Then
$$|p-y|\le |p-q|+|q-v|+|y-v|< 2|p-q|+|y-v|$$ so the open balls centered at
$p,x,y$ with radius equal to the distance to $[u,v]$ leave at most two
segments in $[x,y]$ uncovered, and these segments have length $<|p-q|$
so each can be covered by one more ball.
\end{proof}

\begin{cor}\label{2.3}
Suppose $[x, y]$ and $[u, v]$ are geodesic segments such that
$u\in\pi_{uv}(x)$ and $v\in\pi_{uv}(y)$. Also suppose that $[u,v]$ is
$B$-contracting. Then either $|u-v|<\Phi(B,C)$ or the distance between
$[x,y]$ and $[u,v]$ is $<\Phi(B,C)$.
\end{cor}

\begin{proof}
By Lemma \ref{subsegments contracting} there is
$B'=\Phi_{\ref{subsegments contracting}}(B,C)$ such that
every subsegment of $[u,v]$ is $B'$-contracting. We are done if the
distance $D$ between $[x,y]$ and $[u,v]$ is $\le 5B'+2C$ so suppose it is $>
5B'+2C$. Now pass to a subinterval $[u',v']\subset [u,v]$ of length
$5B'+2C$ whose distance to $[x,y]$ is also $D$. Further, pass
to a subinterval $[u'',v'']\subset [u',v']$ by moving the endpoints at
most $C$ to ensure that $u''\in\pi_{u''v''}(x)$ and
$v''\in\pi_{u''v''}(y)$. The distance between $[u'',v'']$ and $[x,y]$
is $\ge D$ so Lemma \ref{5 balls} applies to $[x,y]$ and $[u'',v'']$
and shows that $|u''-v''|<5B'$, a contradiction.
\end{proof}

\begin{lemma}[Thin Triangles]\label{thin}
Suppose $[a,b]$ is $B$-contracting and $b\in\pi_{ab}(c)$. 
Then $$d(b,[a,c])<3B+C+1$$
\end{lemma}

In the proof we need the following.

\begin{lemma}\label{1}
Suppose $[a,b]$ is $B$-contracting and $b\in\pi_{ab}(c)$. Then 
$$|a-b|+|b-c|\geq |a-c|\geq |a-b|+|b-c|-B-C-1$$
\end{lemma}

\begin{figure}[htbp]
\begin{center}
\input{l1.pstex_t}
\caption{Lemma \ref{thin} and Lemma \ref{1}}
\end{center}
\end{figure}

\begin{proof}
Consider the ball centered at $c$ with radius $|b-c|-1$. By
assumption, the projection of this ball to $[a,b]$ has diameter
$<B$. Let $d$ be the point on $[c,a]$ at distance $|b-c|-1$ from $c$
(if $|b-c|<1$ there is nothing to prove). Denote by $e$ a point in
$\pi_{ab}(d)$. Thus $|b-e|<B$. Now
\begin{equation*}
\begin{split}
|a-c|-|b-c|=|c-d|+|d-a|-|b-c|=|b-c|-1+|d-a|-|b-c|=\\
|d-a|-1\geq
|e-a|-C-1=|a-b|-|b-e|-C-1>|a-b|-B-C-1
\end{split}
\end{equation*}
where the first inequality comes from (DD).
\end{proof}

\begin{proof}[Proof of Lemma \ref{thin}]
Let $q\in [a,c]$ be the point with $|c-q|=|c-b|$. We will argue that
$|b-q|$ is bounded.

Let $r\in \pi_{ab}(q)$. Note that $|b-r|<B$. Subtracting $|b-c|=|c-q|$
from $|a-c|\le |a-b|+|b-c|$ yields $|a-q|\le |a-b|$ and so
$|a-q|<|a-r|+B$. On the other hand, Lemma \ref{1} implies that
$|a-q|>|a-r|+|r-q|-B-C-1$. Comparing these inequalities we deduce
$|r-q|<2B+C+1$. It follows that $|b-q|\le |b-r|+|r-q|<3B+C+1$.
\end{proof}

%For clarity, we will not try to keep up with constants such as
%$3B+C+1$ any more, but rather we introduce the following function
%notation:
%$\Phi=\Phi(B,C)$ stands for a positive function which is monotonically
%increasing in each variable. 
%Sometimes we decorate $\Phi$ with a subscript, and sometimes
%there are more than 2
%arguments. When we want to refer to a $\Phi$-function from a previous
%lemma, we will indicate the number in the subscript,
%e.g. $\Phi_{\ref{thin}}(B,C)=3B+5C+1$ and $\Phi_{\ref{1}}(B,C)=B+C+1$

\begin{lemma}\label{projections}
Let $[a,b]$ be a $B$-contracting segment and assume $|b-c|\le D$. Let
$x\in X$ and denote by $x_1,x_2$ projections of $x$ onto $[a,b],[a,c]$
respectively (i.e. $x_1\in\pi_{ab}(x),x_2\in \pi_{ac}(x)$). Let $x_3$
be a projection of $x_2$ onto $[a,b]$. Then $|x_1-x_3|\le \Phi(B,C,D)$.
\end{lemma}

\begin{figure}[htbp]
\begin{center}
\input{con.pstex_t}
\caption{Lemma \ref{projections}}
\end{center}
\end{figure}

\begin{proof}
Let $x_4$ be a projection of $x_1$ onto $[a,c]$. By Lemma \ref{1}
$$|x_3-x|>|x_3-x_1|+|x_1-x|-\Phi_{\ref{1}}(B,C)$$
By Axiom (FT) we have
$$|x-x_2|\geq |x-x_3|-|x_2-x_3|>|x-x_1|+|x_1-x_3|-\Phi_{\ref{1}}(B,C)-C-D$$
Again by (FT)
$$|x-x_4|\le |x-x_1|+|x_1-x_4|\le |x-x_1|+C+D$$
Since $|x-x_4|\ge |x-x_2|$ we deduce
$$|x-x_1|+C+D>|x-x_1|+|x_1-x_3|-\Phi_{\ref{1}}(B,C)-C-D$$
and the claim follows.
\end{proof}

\begin{lemma}[Stability]\label{stability}
If $[a,b]$ is $B$-contracting, $|a-a'|,|b-b'|\le D$ then $[a',b']$
is $\Phi(B,C,D)$-contracting.
\end{lemma}

\begin{proof}
It suffices to prove this in the special case $a=a'$, for the general
case can be obtained by applying the special case twice. 

Consider a ball disjoint from $[a,b']$. Suppose the projection to
$[a,b']$ has diameter $K$. After shrinking the radius by $C+D$ the
ball is also disjoint from $[a,b]$, its projection $P$ to $[a,b']$ has
diameter $\geq K-(C+D)-2C$, and its projection $Q$ to $[a,b]$ has
diameter $<B$. Every point of $P$ is within $C+D$ of its
projection to $[a,b]$, which in turn is within
$\Phi_{\ref{projections}}(B,C,D)$ of a point in $Q$. Thus $P$ is
contained in the
$C+D+\Phi_{\ref{projections}}(B,C,D)$-neighborhood of $Q$ and also
in any ball of radius $C+D+\Phi_{\ref{projections}}(B,C,D)+B$
centered at a point of $Q$. It follows that $K<
2(C+D+\Phi_{\ref{projections}}(B,C,D)+B)+C+D+2C$.
\end{proof}

\section{Quasi-homomorphisms}\label{s:qh}
In this section we construct many quasi-homomorphisms $\Gamma\to\R$
for a group $\Gamma$ that acts isometrically on the CAT(0) space $X$
satisfying (DD) and (FT).
In Section \ref{wpd}
we will see that many of them are nontrivial (and
linearly independent) under suitable assumptions on the action.

\begin{prop}\label{variation}
Let $[a,b]$ be a $B$-contracting segment and $[p,q]$ an arbitrary
segment. Assume that the distance to $[p,q]$ along $[a,b]$
is minimized at $a$ and is equal to $d\geq 1$. Then
$$d(b,[p,q])-d(a,[p,q])\ge \Big( 1-\frac Bd\Big) | a-b|
-\Phi(B,C).$$
\end{prop}

Thus when $d>>B$ the segment $[a,b]$ is, for the most part,
moving ``orthogonally'' to $[p,q]$.

\begin{proof}
Let $r$ be a projection of $b$ to $[p,q]$, let $s$ be a projection
of $r$ to $[a,b]$ and 
let $t\in [p,q]$ be such that $|a-t|=d$. Subdivide $[t,r]$ into the
smallest possible number, say $K$, of subintervals so that each has
length $<d$. Therefore each subinterval has the projection to $[a,b]$
of diameter $<B$ and we conclude that $|a-s|<KB$. 
We have 
\begin{equation*}
\begin{split}
d(b,[p,q])=|b-r|\ge |b-s|+|s-r|-\Phi_{\ref{1}}(B,C)\ge\\
|b-s|+d-\Phi_{\ref{1}}(B,C)\ge
|a-b|-KB+d-\Phi_{\ref{1}}(B,C)
\end{split}
\end{equation*}
Since $t$ and $r$
are projections of $a$ and $b$, (DD) implies that
$$|a-b|\ge |t-r|-C\ge (K-1)d-C$$
and we deduce that $K\le 1+\frac{|a-b|+C}d$. Therefore
$$d(b,[p,q])-d>|a-b|-KB-\Phi_{\ref{1}}(B,C)\ge
|a-b|-\Big(1+\frac{|a-b|+C}d\Big)B-\Phi_{\ref{1}}(B,C)$$ and the claim
follows.
\end{proof}

\subsection{Expressways}

Assume a group $\Gamma$ acts on $X$ by isometries.  Fix a collection
$\mathcal E$ of geodesic segments in $X$ (thought of as {\it
  expressways}), subject to the following axioms.  The constants $L$
and $B$ are fixed ($L$ will be chosen large using Proposition
\ref{Hausdorff neighborhoods}).
\begin{itemize}
\item $\mathcal E$ is $\Gamma$-equivariant,
\item each $\sigma\in \mathcal E$ is $B$-contracting,
\item each $\sigma\in\mathcal E$ has length $L$.
\end{itemize}

To each piecewise geodesic path $\alpha=[a_0,a_1]\cup [a_1,a_2]\cup\cdots\cup
  [a_{k-1},a_k]$ in $X$ we will assign {\it modified
  length}, equal to $ml(\alpha)=\sum_{i=1}^k |a_{i-1}-a_i|-e$ where $e$ is the
  number of segments $[a_{i-1},a_i]$ which are expressways.
  We can think of modified length as time needed to
  travel along the path, provided we travel with speed 1 on regular
  roads and with fixed speed $>1$ on expressways.

We will call a piecewise geodesic path {\it admissible} if for any two
consecutive geodesic segments $[a_{i-1},a_i]\cup [a_i,a_{i+1}]$, at
least one is an expressway. Note that if neither is an expressway,
then replacing $[a_{i-1},a_i]\cup [a_i,a_{i+1}]$ by
$[a_{i-1},a_{i+1}]$ does not increase modified length, and a repeated
application of this operation converts the path to an admissible path.

Now fix a (long) geodesic segment $[p,q]$. In the next
Proposition we will study admissible paths that join $p$ to $q$ and
minimize modified length (with respect to a fixed collection $\mathcal
E$ of expressways). When $X$ is a proper metric space it is easy to
see that such minimizers exist; in general we study admissible paths
whose modified length is close to minimal. We will show that such
paths are necessarily contained in a Hausdorff neighborhood of
$[p,q]$.

\begin{prop} \label{Hausdorff neighborhoods}
There is a function $\Phi(B,C)$ with the
  following property.  Assume that expressways have length
$L>\Phi(B,C)$. Then for every segment $[p,q]$ every admissible path that
  minimizes modified distance
from $p$ to $q$ is contained in the Hausdorff
$\Phi(B,C)$-neighborhood of $[p,q]$. The same is true (even if no
minimizing paths exist) for each admissible path within 1 of the infimum.
\end{prop}

In the proof we need the following lemma.

\begin{lemma}\label{needed for Hausdorff neighborhoods}
Let $[t,t']$ be a geodesic segment in $X$ and let $[t,y]\cup [y,z]\cup
[z,t']$ be a piecewise geodesic joining $t$ and $t'$. Also assume that
every subsegment of $[t,y]$ and of $[z,t']$ is $B'$-contracting. Then
either
\begin{enumerate}[(1)]
\item $|t-y|+|y-z|+|z-t'|-|t-t'|>3$, or
\item $[t,y]\cup [y,z]\cup
[z,t']$ is contained in the $\Phi(B',C)$-neighborhood of $[t,t']$.
\end{enumerate}
\end{lemma}

\begin{proof}
We may assume that one of $[t,y]$ or $[z,t']$, say $[z,t']$, is
degenerate (a point), since the general case follows by applying the
special case twice.

Put $T=\Phi_{\ref{thin}}(B',C)$.
Let $\tilde t \in \pi_{ty}(t')$,
$t_1 \in \pi_{tt'}(\tilde t)$ and 
$t_2 \in \pi_{yt'}(\tilde t)$.
By Lemma \ref{thin},
$|\tilde t - t_1| < T, |\tilde t -t_2 | < T$.

 The ``savings'' $|t-y|+|y-t'|-|t-t'|$
is equal to 
$(|t-\tilde t| -|t-t_1|)
+(|t_2-t'|-|t_1-t'|)
+(|\tilde t -y|+|y-t_2|)$,
which is at least
$-T-2T+(2|\tilde t-y| -T)$.
Therefore, if (1) does not hold, then this 
number is less than $3$, so that 
$|\tilde t -y| < 2T +\frac{3}{2}$.
It follows that 
$[t,y]$ is within 
distance
$T+C+(2T+\frac{3}{2})$ from $[t,t']$.
Also, $[y,t']$ is 
within distance
$2T+C+(2T+\frac{3}{2})+T$ from 
$[t,t']$.
Put $\Phi(B',C)=5T+C+\frac{3}{2}$, then 
(2) holds.
\end{proof}

\begin{proof}[Proof of Proposition \ref{Hausdorff neighborhoods}]
Fix a (nearly) minimizing admissible path $\alpha=[a_0,a_1]\cup
[a_1,a_2]\cup\cdots\cup [a_{k-1},a_k]$ from $p=a_0$ to $q=a_k$. For
each expressway $\sigma=[a_i,a_{i+1}]$ in $\alpha$ let
$\Delta(\sigma)$ denote the distance between $\sigma$ and
$[p,q]$. Also, let
$B'=\Phi_{\ref{subsegments contracting}}(B,C)$ so that subsegments of
expressways in $\alpha$ are $B'$-contracting.

There are now two cases.

\noindent
{\bf Case 1.} Maximal $\Delta(\sigma)$ is $\ge \max(1,2B')$. Say the
maximum is realized on the expressway $\sigma=[a_i,a_{i+1}]$ and let
$t\in \sigma$, $\hat t\in [p,q]$ be such that $|t-\hat
t|=\Delta(\sigma)$. For concreteness, say $|t-a_{i+1}|\ge
|t-a_i|$. Then consider the subpath $[a_i,a_{i+1}]\cup
[a_{i+1},a_{i+2}]\cup [a_{i+2},a_{i+3}]=[x,y]\cup [y,z]\cup [z,w]$ (we
certainly cannot have $i+1=k$; if $i+2=k$ then
$[a_{i+2},a_{i+3}]=[z,w]$ is degenerate). 
Let $t'\in [z,w]$ be a point
closest to $[p,q]$. By the choice of $\sigma$, $d(t,[p,q])\ge
d(t',[p,q])$ and by Axiom (FT) $[t,t']$ is contained in the
$C+\Delta(\sigma)$-neighborhood of $[p,q]$. Now we may apply Lemma
\ref{needed for Hausdorff neighborhoods}. If (1) holds, then replacing
the portion $[x,y]\cup [y,z]\cup [z,w]$ of $\alpha$ by $[x,t]\cup
[t,t']\cup [t',w]$ results in a path $\alpha'$ with modified length
$ml(\alpha')< ml(\alpha)-1$. The path $\alpha'$ may not be admissible,
but we may replace it by an admissible path without increasing
modified length. Thus $\alpha$ wasn't within 1 of the infimum. Therefore,
(2) must hold and in particular $y$ is in the
$d(t,[p,q])+C+\Phi_{\ref{needed for Hausdorff
neighborhoods}}(B',C)$-neighborhood of $[p,q]$. It follows that the
distance function to $[p,q]$ has variation $\le C+\Phi_{\ref{needed for
Hausdorff neighborhoods}}(B',C)$ along $[t,y]$. But this variation is
at least $\frac 12 |t-y|-\Phi_{\ref{variation}}(B',C)$ and we conclude
$$L=|x-y|\le 2|t-y|\le 4(\Phi_{\ref{variation}}(B',C)+C+\Phi_{\ref{needed for
Hausdorff neighborhoods}}(B',C))$$

\noindent
{\bf Case 2.} Every expressway in $\alpha$ intersects the Hausdorff
$\max(1,2B')$-neighborhood of $[p,q]$. For each such expressway
$[a_i,a_{i+1}]$ let $[s_i,t_i]$ be the maximal subinterval with
$d(s_i,[p,q])\le \max(1,2B')$ and $d(t_i,[p,q])\le \max(1,2B')$. Let
$D_i$ be the length of the longer of the two complementary components
(one or both of these might be empty; the length of the empty set is
0). Arguing exactly as in Case 1, with $[t,y]=[t_i,a_{i+1}]$ or
$[t,y]=[s_i,a_i]$, we conclude
$$D_i\le 2(\Phi_{\ref{variation}}(B',C)+\Phi_{\ref{needed for
Hausdorff neighborhoods}}(B',C))$$ (or else $\alpha$ is not within 1
of infimum). Therefore, $\alpha$ is in the Hausdorff 
$2(\Phi_{\ref{variation}}(B',C)+\Phi_{\ref{needed for
Hausdorff neighborhoods}}(B',C))+\max(1,2B')$-neighborhood of $[p,q]$.
\end{proof}

For the rest of this section we will use the following notation.
Let $D=\Phi_{\ref{Hausdorff neighborhoods}}(B,C)$, 
$S=\Phi_{\ref{stability}}(B,C,D)$, $S'=\Phi_{\ref{subsegments
    contracting}}(S,C)$, and $T=\Phi_{\ref{2.3}}(S',C)+D$.

\begin{lemma}\label{no expressways in thick part}
Every $B$-contracting segment $[x,y]$ of length $\ge\Phi(B,C)$ contained in
the $D$-neighborhood of $[p,q]$ contains a point whose 
distance from $[p,r]\cup [r,q]$ is $<T$ for any $r \in X$.
\end{lemma}

\begin{proof} 
Let $u$ be a projection of $x$ and $v$ a projection of $y$ to
$[p,q]$. Then $[u,v]$ is $S$-contracting and every subsegment is
$S'$-contracting. Say a
projection of $r$ to $[u,v]$ is closer to $v$ than to $u$.
Assuming $|x-y|$ is sufficiently big, it follows 
that some point of $[u,v]$ is within $\Phi_{\ref{2.3}}(S',C)$ of $[p,r]$.
\end{proof}

\subsection{Construction of quasi-homomorphisms}

Now assume that a group $\Gamma$ acts on $X$ by isometries.
Fix a constant $L$ (in the 
definition of $\mathcal E$) to satisfy 
the assumption in Proposition \ref{Hausdorff neighborhoods}.
Let $\mathcal E$ be a collection of expressways. For $a, b\in  X$ denote by
$$\lambda(a, b) = \lambda_{\mathcal E}(a, b)$$ 
the infimum of modified lengths of admissible piecewise
geodesic paths from $a$ to $b$.

\begin{prop}\label{lambda}
\begin{itemize}
\item $\lambda(a,b)\leq |a-b|$.
\item $|\lambda(a,b)-\lambda(a',b')|\leq |a-a'|+|b-b'|$.
\item $\lambda(ga,gb)=\lambda(a,b)$ for any $g\in \Gamma$.
\item If there are no expressways in the $D$-neighborhood of
  $[a,b]$ then $\lambda(a,b)=|a-b|$.
\item $\lambda(a,c)\leq\lambda(a,b)+\lambda(b,c)$ and if $b\in [a,c]$
  then also $\lambda(a,c)>\lambda(a,b)+\lambda(b,c)-2D-1$.
\end{itemize}
\end{prop}

\begin{proof}
The first three statements follow immediately from the definition, the
fourth one follows from Proposition \ref{Hausdorff neighborhoods}. The
first part of the last one is obvious, while the second part also
follows from Proposition \ref{Hausdorff neighborhoods} by breaking up
the (nearly) minimizing path from $a$ to $c$ at a point near $b$.
\end{proof}

Now fix a basepoint $x_0\in X$ and a $B$-contracting segment
$\sigma$. 
Let $\mathcal{E}$ consist of all translates of $\sigma$ by $\Gamma$.
Define
$$\phi_\sigma:\Gamma\to\R$$
via
$$\phi_\sigma(g)=\lambda(g(x_0),x_0)-\lambda(x_0,g(x_0))$$    

\begin{thm}\label{Bcont.qh}
If $|\sigma|>\Phi(B,C)$ then $\phi_\sigma$ is a quasi-homomorphism,
i.e. $$|\phi_\sigma(gg')-\phi_\sigma(g)-\phi_\sigma(g')|$$ is uniformly bounded
over all $g,g'\in \Gamma$ (by a function of $B$ and $C$).
\end{thm}

\begin{proof}
Consider the triangle with vertices $x_0,g(x_0),gg'(x_0)$. We have
\[
\begin{array}{ll}
\phi_\sigma(g)&=\lambda(g(x_0),x_0)-\lambda(x_0,g(x_0))\\
\phi_\sigma(g')&=\lambda(gg'(x_0),g'(x_0))-\lambda(g'(x_0),gg'(x_0))\\
\phi_\sigma(gg')&=\lambda(gg'(x_0),x_0)-\lambda(x_0,gg'(x_0))
\end{array}
\]
Now consider the following ``thick-thin'' decomposition of the
triangle. Each side is subdivided into 3 subintervals (the case when
one or more middle intervals is degenerate is easier and is left to
the reader). Figure \ref{thick-thin} indicates the naming
of the subdivision points. The point $a_1$ is chosen to be the last
point in $[x_0,g(x_0)]$ whose distance to $[x_0,gg'(x_0)]$ is $\le
T+D$, and similarly for the other 5 subdivision points.

Note that $|a_1-a_2|\le C+T+D+2T+2D$ and similarly
for the other two pairs. (To see this, let $t$ be the projection of
$a_2$ onto $[x_0,g(x_0)]$. Then $t\in [x_0,a_1]$ by the choice of
$a_1$. We will argue that $|t-a_1|\le C+T+D+T+D$ and hence $|a_1-a_2|\le
C+T+D+2T+2D$. Let $s$ be
a projection of $a_1$ onto $[x_0,gg'(x_0)]$ (hence $s\in
[x_0,a_2]$). By Axiom (FT) applied to $[x_0,a_2]$ and $[x_0,t]$ we
see that $s$ is at distance $\le C+T+D$ from some point $r\in
[x_0,t]$. Thus $|r-a_1|\le |r-s|+|s-a_1|\le C+T+D+T+D$ and since
$t\in [r,a_1]$ the claim is proved.)

\begin{figure}[htbp]
\begin{center}
\input{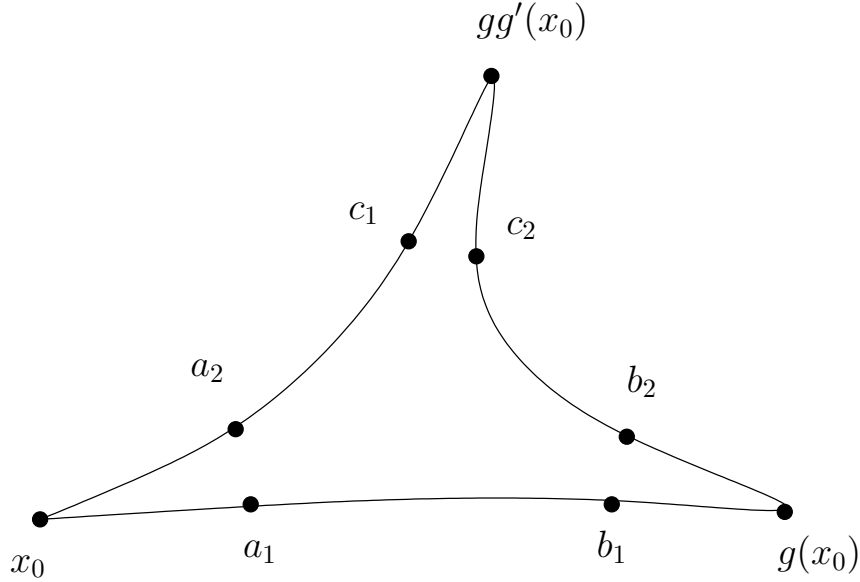}
\caption{Thick-thin decomposition}
\end{center}
\label{thick-thin}
\end{figure}

Now replace each of the 6 terms above by the sum of 3 terms using the
subdivision, e.g.
$$\lambda(g(x_0),x_0)$$ by the sum 
$$\lambda(g(x_0),b_1)+\lambda(b_1,a_1)+\lambda(a_1,x_0).$$ Each such
replacement introduces an error that does not exceed $4D+2$ by
Proposition \ref{lambda}.  Then consider the sum
$\phi_\sigma(gg')-\phi_\sigma(g)-\phi_\sigma(g')$. The terms such as
$\lambda(a_1,b_1)-\lambda(b_1,a_1)$ cancel in pairs since they equal
the distance $|a_1-b_1|$ (by Proposition \ref{lambda} and Lemma
\ref{no expressways in thick part}). 
To see it, assume that $[x,y]$ is an expressway contained
in the $D$-neighborhood of $[a_1,b_1]$.
By Lemma \ref{no expressways in thick part} 
with $x_0=p,g(x_0)=q$ and $gg'(x_0)=r$,
there is a point $s$ in $[x,y]$ at distance 
$<T$ of $[x_0,gg'(x_0)]$ or $[g(x_0),gg'(x_0)]$,
say $[x_0,gg'(x_0)]$.
There is a point $s'$ in $[a_1,b_1]$ within $D$ of $s$.
But then $s'$ is at distance $<D+T$ of $[x_0,gg'(x_0)]$, 
contradicting the choice of $a_1$.
The terms such as
$\lambda(x_0,a_1)-\lambda(x_0,a_2)$ also cancel, up to a bounded error
by the above remark. At the end, all terms cancel and the
total error is bounded as a function of $B$ and $C$.
\end{proof}

\section{Rank 1 isometries}\label{rk1}

We say that two subsets of a CAT(0) space $X$ are {\it $B$-Hausdorff
equivalent} if each is contained in the $B$-Hausdorff neighborhood of
the other. We continue to assume that $X$ satisfies Axioms (DD) and
(FT). 

\begin{definition}
Let $g:X\to X$ be an isometry and let $x_0\in X$ be a
basepoint. Denote $x_n=g^n(x_0)$. We say that $g$ is {\it rank
1} if $|x_n-x_0|\to\infty$ as $n\to\infty$
and there exists $B>0$ such that for every $n>0$
\begin{itemize}
\item some geodesic $[x_0,x_n]$ is 
  $B$-Hausdorff equivalent to $\{x_0,x_1,\cdots,x_n\}$, and
\item some geodesic $[x_0,x_n]$ is $B$-contracting.
\end{itemize}
\end{definition}

Using Axiom (FT) and Lemma \ref{stability} it is easy to see that
the notion does not depend on the choice of the basepoint or the
geodesics. It is also easy to see that there is $\epsilon>0$ such that
$|x_0-x_n|\ge n\epsilon$ for all $n>0$.

\begin{prop}
If $X$ is a $\delta$-hyperbolic space then every hyperbolic isometry
$X\to X$ is rank 1.\qed
\end{prop}

\subsection{Characterization when $X$ is a proper $CAT(0)$ space}
Ballmann and Brin \cite{bb} defined rank 1 isometries in the case when
$X$ is a proper $CAT(0)$ space, namely, a hyperbolic isometry $g:X\to
X$ is rank 1 if some (equivalently every) axis of $g$ fails to bound a
half-flat.

To reconcile the two definitions, we show they are equivalent in this
situation. 
We will need the following statement from the work of Ballmann, see
the proof of Lemma III.3.3
in \cite{bbook}.
By $\overline X$ denote the usual visual compactification of
$X$.

\begin{prop}\label{bbinput}
Let $\ell$ be an axis of $g$ and assume that $g$ is rank 1. Let $\{x_n\}_n$
be a sequence of points in $X$ such that the projections of $x_n$ to
$\ell$ leave every compact set. Then the accumulation set of the
sequence $\{x_n\}$ in $\overline X$ is contained in the set of two
endpoints of $\ell$.\qed
\end{prop}
 
\begin{thm}
Let $X$ be a proper $CAT(0)$ space. Let $g:X\to X$ be a hyperbolic
isometry with axis $\ell$. Then $\ell$ is $B$-contracting for some $B$
if and only if $\ell$ fails to bound a half-flat.
\end{thm}

\begin{proof}
If $\ell$ bounds a half-flat then $\ell$
is clearly not $B$-contracting for any $B$. 

Now assume that $\ell$ is not $B$-contracting for any $B$. Then there
is a sequence of balls in the complement of $\ell$ whose projections
to $\ell$ are larger and larger intervals. Denote by $\{z_n\}$ the
sequence of centers of the balls. We may assume that the projection of
$z_n$ to $\ell$ is contained in a fixed closed interval $J$
(a fundamental domain for $g:\ell\to\ell$). Let $x_n$ be a point in the
$n$th ball chosen so that the sequence of projections of the $x_n$'s
to $\ell$ leave every compact set. By Proposition \ref{bbinput}, after
passing to a subsequence we may assume that the sequence $\{x_n\}$
converges to an endpoint $t\in\overline X$ of $\ell$, and also that
the sequence $\{z_n\}$ converges to a point $\overline z\in\overline
X$.

Moreover, we may assume that there exists some constant $r>0$ such that 
$d(x_n,\ell) \le r$ for all $n$.
This is because, otherwise, we can retake $x_n$
to satisfy this. 
Indeed, since $\ell$ does not bound a flat half plane, it does
not bound a flat strip of width, say, $r>0$, since $X$ is proper.
Let $Z_n \in J$ be the projection of $z_n$ to $\ell$.
Let $x_n'$ be the first point of the geodesic 
from $x_n$ to $Z_n$ to be contained in the $r$-neighborhood
of $\ell$. (If $x_n$ is contained in the $r$-neighborhood,
then $x_n'=x_n$.)
Then by \cite[Lemma 3.1]{bbook} combined with 
\cite[Lemma 3.3]{bbook}, the projections of the $x_n'$'s
to $\ell$ leave every compact set.
We keep denoting those points by $x_n$ in the rest.

Consider the sequence of functions $f_n:X\to\R$ given by
$$f_n(y)=|z_n-y|-d(z_n,\ell)$$
This sequence converges to a Busemann function $f:X\to\R$ centered at
$\overline z$. The horoball $f^{-1}(-\infty,0]$ is convex, it includes a
point of $J\subset \ell$, and its closure in $\overline X$ includes
$t$. It follows that the ray $R$ of $\ell$ joining a point of $J$ with $t$
is contained in this horoball. 
(Indeed, $f_n \le 0$ on $[x_n,Z_n]$. Since $x_n$ is contained
in the $r$-neighborhood of $\ell$ and $Z_n \in J$,
the sequence $[x_n,Z_n]$ converges to $R$.
Thus $f \le 0$ on $R$.) 
Moreover, each function $f_n$ is
1-Lipschitz and
non-negative on $\ell$, so the same holds for $f$. In particular, $f$
is identically 0 on the ray $R$.

Let $p\in R$ and consider the ray $l_p$ joining $p$ with $\overline
z$. If we isometrically identify $l_p=[0,\infty)$ then $f(y)=-y$ for
  $y\in l_p$. It follows that the projection of $l_p$ onto $\ell$ is
  the point $p$. 

Now take $p,q\in R$ and consider the function $g:l_p=[0,\infty)\to\R$
  given by
$$g(y)=d(y,l_q)$$
This function is non-increasing because this is true for the analogous
function defined with respect to two radial line segments sharing
their terminal endpoints. We also have that $g(y)\geq |p-q|$ since
the projection to $\ell$ doesn't increase distances and $l_p$ and
$l_q$ project to $p$, $q$ respectively. It follows that $g$ is a
constant function and hence $l_p$ and $l_q$ cobound a flat half-strip.
In this way we can construct arbitrarily wide
half-strips by taking $|p-q|$ large.

Now using the properness of $X$, a standard limiting argument applied to
the translates of these half-strips by the powers of $g$ produces a
half-flat bounded by $\ell$.
\end{proof}

\subsection{Schottky groups}

\begin{lemma}\label{no crossing}
Suppose $|q-q'|\le D$, $s\in [x,q]$, $s'\in [x,q']$. Then either
$d(s,[s',q'])\le \Phi(C,D)$ or $d(s',[s,q])\le \Phi(C,D)$.
\end{lemma}

\begin{proof}
Let $r$ be a projection of $s'$ to $[x,q]$ and $r'$ a projection of
$s$ to $[x,q']$. If $r'\in [s',q']$ or if $r\in [s,q]$ we are done by
Axiom (FT) ($|r-s'|$ and $|r'-s|$ are bounded), so let us
suppose that $r'\in [x,s']$ and $r\in [x,s]$. Now apply 
(FT) to $[x,r]$ and $[x,s']$ to deduce that $d(r',[x,r])$
is bounded. Thus $|r-s|$ is bounded and so is $|s-s'|$.
\end{proof}

\begin{lemma}\label{projections in an interval}
Let $\sigma$, $\sigma'$ be two $B$-contracting segments. Suppose that
the minimal distance between $\sigma$ and $\sigma'$ is $D$. Let $x\in
X$ be an arbitrary point and let $p,p'$ be projections of $x$ to
$\sigma,\sigma'$ respectively. Then either $d(p,\sigma')\le
\Phi(B,C,D)$ or $d(p',\sigma)\le
\Phi(B,C,D)$.
\end{lemma}

\begin{figure}[htbp]
\begin{center}
\input{schottky3.pstex_t}
\caption{Lemma \ref{no crossing} and Lemma \ref{projections in an interval}}
\end{center}
\label{schottky3}
\end{figure}

\begin{proof}
Let $q\in\sigma$, $q'\in\sigma'$ be points with $|q-q'|=D$. Let $s$ be
a projection of $p$ to $[x,q]$ and $s'$ a projection of $p'$ to
$[x,q']$. By Lemma \ref{thin} both $|p-s|$ and $|p'-s'|$ are bounded.
By Lemma \ref{no crossing} we may assume, by symmetry,
that $d(s,[s',q'])$ is bounded. Say $\hat s\in [s',q']$
is such that $|s-\hat s|$ is bounded. By Axiom (FT) $d(\hat
s,[p',q'])$ is bounded, so $d(p,\sigma')$ is bounded.
\end{proof}

A slight sharpening, in case $\sigma$ and $\sigma'$ share an endpoint,
is the following:

\begin{lemma}\label{sharpening}
Suppose $[x,y]$ and $[y,z]$ are $B$-contracting segments and $p\in
X$. Let $a$ be a projection of $p$ to $[x,y]$ and $b$ a projection of
$p$ to $[y,z]$. Then either:
\begin{itemize}
\item $d(a,[y,z])\le \Phi(B,C)$ and $|p-a|\ge |p-b|-\Phi(B,C)$, or
\item $d(b,[x,y])\le \Phi(B,C)$ and $|p-b|\ge |p-a|-\Phi(B,C)$.
\end{itemize}
\end{lemma}

\begin{figure}[htbp]
\begin{center}
\input{sharpening.pstex_t}
\caption{Lemma \ref{sharpening}}
\end{center}
\label{schottky3-fig}
\end{figure}

\begin{proof}
Let $s$ be a projection of $a$ onto $[p,y]$ and $t$ a projection of
$b$ onto $[p,y]$. Then either $t\in [p,s]$ or $s\in [p,t]$, say the
former. Note that both $|a-s|$ and $|b-t|$ are bounded 
by Lemma \ref{thin}. Consider the segments $[y,t]$ and $[y,b]$. Since
$|b-t|$ is bounded, Axiom (FT) implies that $d(s,[y,b])$ is
bounded as a function of $B$ and $C$, so we deduce $d(a,[y,z])$
is bounded. In addition, $|p-b|\le |p-t|+|b-t|\le |p-s|+|b-t|\le
|p-a|+|a-s|+|b-t|$ so the first bullet holds.
\end{proof}

\begin{definition}\label{independent}
Let $g,h:X\to X$ be two rank 1 isometries. Denote
$x_n=g^n(x_0)$ and $y_n=h^n(x_0)$. We say that
$g$ and $h$ are {\it independent} if the function $\Z^2\to
[0,\infty)$ given by $(m,n)\mapsto |x_m-y_n|$ is proper.
\end{definition}

\begin{prop}\label{ping-pong}
Let $g,h$ be two independent rank 1 isometries. Then for
sufficiently large even $N$ the powers $g^N,h^N$ generate a
nonabelian free subgroup of the isometry group of $X$.

Moreover, one can arrange that for any given $E>0$ and for every
reduced word $w$ of length $|w|$ in $g^N$ and $h^N$ we have
$$|x_0-w(g^N,h^N)(x_0)|\ge |w|E$$
\end{prop}

\begin{proof}
We retain the notation from Definition \ref{independent}.  Let
$I_k=[x_{-k},x_k]$ and $J_k=[y_{-k},y_k]$. When $k$ is large, the
shortest distance between $I_k$ and $J_k$ is realized by a pair of
points each of which is near the middle of its respective
interval. From Lemma \ref{projections in an interval} it follows that
(for $k$ very large) every $x\in X$ has projections to either $I_k$ or
to $J_k$ a bounded distance away from the middle of the interval; in
any case, we can ensure that it is contained in the middle third $\hat
I_k$ or $\hat J_k$ of the interval $I_k$ or $J_k$ by making $k$
large. Using the fact that $|x_0-x_n|\to\infty$ we can also ensure
that the distance between $\hat I_k$ and $g^{\pm k}(\hat I_k)$ is large,
and similarly, the distance between $\hat J_k$ and $h^{\pm k}(\hat J_k)$
is large. Let $U_{\pm}$ ($V_{\pm}$) be the subset of $X$ consisting of
points that have a projection to $I_k$ ($J_k$) contained in the
component of $I_k-\hat I_k$ ($J_k-\hat J_k$) that contains $x_{\pm k}$
($y_{\pm k}$). Note that $U_-\cap U_+=\emptyset$ as long as $\hat I_k$
is sufficiently long (which is ensured by making $k$ large), and
likewise $V_-\cap V_+=\emptyset$. From Lemma \ref{projections in an
  interval} it follows that $U_{\pm}\cap V_{\pm}=\emptyset$.

{\bf Claim:} $g^{2k}(X-U_-)\subset U_+$. Indeed, suppose not. Then
some $x\in X$ has a projection to $I_k=[x_{-k},x_k]$ contained in the
middle or the last third of $I_k$ while $g^{2k}(x)$ has its
projection to $I_k$ contained in the first two thirds. Note that
$g^{2k}(x)$ has a projection to $g^{2k}(I_k)$ contained in the
last two thirds of the interval, so we have a contradiction to Lemma
\ref{projections in an interval} (for $k$ large). 

In a similar way one verifies that $g^{-2k}(X-U_+)\subset U_-$ and 
$h^{{\pm}2k}(X-V_{\mp})\subset V_{\pm}$. Thus, by the standard
ping-pong argument, $<g^{2k},h^{2k}>$ is free.

For the ``moreover'' part, arrange that the $E$-neighborhood of
$g^{{\pm}2k}(X-U_{\mp})$ is contained in $U_{\pm}$ and that the
$E$-neighborhood of
$h^{{\pm}2k}(X-U_{\mp})$ is contained in $U_{\pm}$. 
\end{proof}

\begin{lemma}\label{qg}
There is a function $\Phi(B,C)$ with the following property.
Let $\gamma=[x_0,x_1]\cup [x_1,x_2]\cup\cdots\cup [x_{k-1},x_k]$ be a
piecewise geodesic such that each $[x_i,x_{i+1}]$ is
$B$-contracting. Assume also that
$$d([x_{i},x_{i+1}],[x_{i+2},x_{i+3}])>\Phi(B,C)$$ for
$i=0,1,\cdots,k-3$. Then 
\begin{itemize}
\item $[x_0,x_k]$ is contained in the
$\Phi(B,C)$-neighborhood of $\gamma$.
\item
$[x_0,x_k]$ is $\Phi(B,C)$-contracting.
\end{itemize}
\end{lemma}

\begin{proof}
Consider the projections of points in $[x_0,x_k]$ onto
$[x_{i-1},x_{i}]$ and onto $[x_{i},x_{i+1}]$. At least one of the two
sets is contained in the $\Phi_{\ref{sharpening}}(B,C)$-neighborhood
of the other geodesic. Since by Axiom (DD) the projection is coarsely
continuous, we see that there is at least one point $t\in [x_0,x_k]$
such that both projections to $[x_{i-1},x_{i}]$ and to
$[x_{i},x_{i+1}]$ are in the
$\Phi_{\ref{sharpening}}(B,C)+C$-neighborhood of the other geodesic
(when $|x_j-x_{j+1}|$ are all large, the projection of $x_0$ to each
$[x_j,x_{j+1}]$ is close to $[x_{j-1},x_j]$ by induction on $j$ using
Lemma \ref{sharpening}, and similarly the projection of $x_k$ to each
$[x_j,x_{j+1}]$ is close to $[x_{j+1},x_{j+2}]$). Denote by $s_i$ the
last point along $[x_0,x_k]$ with that property. Also let $s_0=x_0$
and $s_k=x_k$ so that $s_0<s_1<\cdots<s_k$.

Now fix some $N=N(B,C)$ and consider the $N$-neighborhood $\mathcal N$
of $\gamma$. Suppose $[t,t']$ is a component of $[x_0,x_k]-\mathcal
N$. For concreteness, say $t\in [s_{i-1},s_i]$ and $t'\in
[s_{j-1},s_j]$. We see immediately that if $N(B,C)$
is sufficiently large, $|j-i|\le 1$
since otherwise Corollary \ref{2.3} is violated (for the projection of
$[s_i,s_{i+1}]$ to $[x_i,x_{i+1}]$). Say $j=i+1$ for
concreteness (the case $j=i$ is easier). The diameters
of the projections of $[t,s_i]$ and of $[s_i,t']$ to $[x_{i-1},x_i]$
and $[x_i,x_{i+1}]$ are bounded again by Corollary \ref{2.3}, so
$|t-t'|$ is bounded by $2N$ plus the diameters of the projections plus
the jump from one projection to the other. We conclude that every
subinterval of $[x_0,x_k]$ which is outside $\mathcal N$ has bounded
length, so $[x_0,x_k]$ is in a bounded neighborhood of $\gamma$.

By Lemma \ref{stability} each $[s_i,s_{i+1}]$ is
$B'$-contracting for some $B'=B'(B,C)$. By choosing $F(B,C)$ large we
can arrange that $|s_i-s_{i+1}|$ are guaranteed to be much bigger than
$2B'$. It now follows that $[x_0,x_k]$ is $2B'+C$-contracting (if the
projection of some ball missing $[x_0,x_k]$ to $[x_0,x_k]$ has
diameter $\ge 2B'+C$ then the projection of the same ball to some
$[s_i,s_{i+1}]$ has diameter $\ge B'$ since the gaps in the projection
are bounded by $C$). 
\end{proof}

Proposition \ref{ping-pong} and Lemma \ref{qg} immediately give:

\begin{prop}\label{schottky}
Let $g,h$ be two independent rank 1 isometries of $X$. Then
there is $N>0$ such that the group generated by $g^N$ and $h^N$
is nonabelian and free, and every nontrivial element is rank
1. Moreover, there is $B'>0$ such that $[x_0,w(g^N,h^N)(x_0)]$
is $B'$-contracting for every reduced word $w$.\qed
\end{prop}

\section{Weak Proper Discontinuity and Infinite Dimension of $\widetilde
  {QH}(G)$}\label{wpd}

This section is a straightforward adaptation of \cite{bf}. Assume
$\Gamma$ acts on a CAT(0) space $X$ by isometries, and $X$ satisfies
Axioms (DD) and (FT). As before, $x_0\in X$ is a fixed basepoint.

\begin{definition}
For a pair of rank 1 elements $g,h\in \Gamma$ 
write $g\sim h$ if there
is $K>0$ and sequences $m_i,n_i\to\infty$ and $\gamma_i\in\Gamma$
such that $[x_0,g^{m_i}(x_0)]$ and $\gamma_i[x_0,h^{n_i}(x_0)]$
are $K$-Hausdorff equivalent.
\end{definition}

It is easy to see that this is an equivalence relation and that the
concept does not depend on the choice of $x_0$.

\begin{prop}\label{p2}
Suppose there exist
independent rank 1 elements $g_1,g_2\in \Gamma$ such that $g_1\not\sim
g_2$. Then there is a sequence $f_1,f_2,\cdots\in \Gamma$ of rank 1
elements such that
\begin{itemize}
\item $f_i\not\sim f_i^{-1}$ for $i=1,2,\cdots$, and
\item $f_i\not\sim f_j^{\pm 1}$ for $j<i$.
\end{itemize}
\end{prop}

\begin{proof}
Given Proposition \ref{schottky}, the proof is identical to the proof
of Proposition 2 in \cite{bf}.
\end{proof}

\begin{thm}\label{t1}
Suppose there exist
independent rank 1 elements $g_1,g_2\in \Gamma$ such that $g_1\not\sim
g_2$. Then $\widetilde{QH}(\Gamma)$ is infinite-dimensional.
\end{thm}

\begin{proof}
Let $f_1,f_2,\cdots$ be as in Proposition \ref{p2}. We may assume that
each $f_i$ is cyclically reduced as a word in $g_1^N$ and $g_2^N$
($N$ is the constant from Proposition \ref{schottky}). 
Choose
a sufficiently rapidly growing sequence $a_i\in\Z$ and let
$\phi_i:\Gamma\to\R$ be the quasi-homomorphism $\phi_\sigma$ from Theorem
\ref{Bcont.qh} associated to a segment $\sigma=[x_0,f_i^{a_i}(x_0)]$.
The proof that $\phi_1,\phi_2,\cdots$ are linearly independent in
$\widetilde{QH}(G)$ is identical to the proof of Theorem 1 of
\cite{bf}.

A quick summary of the argument is that 
$\phi_i$ is a quasi-homomorphism since $f_i$ is rank 1
(Theorem \ref{Bcont.qh}), 
and it is unbounded on $\langle f_i \rangle$ because $f_i\not\sim f_i^{-1}$.
Indeed (see also \cite{bf}), we can easily arrange that 
$f_i \in [\G,\G]$ in Proposition \ref{p2}.
It then immediately follows that $\phi_i$ is 
non-trivial in $\widetilde{QH}(\G)$.
Also, $\phi_j$ is bounded
on $<f_i>$ for $j<i$ by $f_i\not\sim f_j^{\pm 1}$ which implies
linear independence
in $\widetilde{QH}(\G)$.
\end{proof}

\begin{definition}\label{def.wpd}
We say that the action of $\Gamma$ on $X$ satisfies $WPD$ (weak proper
discontinuity) if
\begin{itemize}
\item $\Gamma$ is not virtually cyclic,
\item $\Gamma$ contains at least one rank 1 element, and
\item for every rank 1 element $g\in \Gamma$, and every
  $c>0$ there exists $M>0$ such that the set
$$\{\gamma\in \Gamma||x_0-\gamma(x_0)|\leq c,|g^M(x_0)-\gamma
g^M(x_0)|\le c\}$$ is finite.
\end{itemize}
\end{definition}

The concept doesn't depend on the choice of the basepoint.

\begin{prop}\label{p6}
Suppose the action of $\Gamma$ on $X$ satisfies $WPD$. Then
\begin{enumerate}[(1)]
\item for every rank 1 element $g\in \Gamma$ the centralizer $C(g)$ is
  virtually cyclic,
\item the action of $\Gamma$ on $X$ is nonelementary (i.e. no line is
  preserved by the whole group),
\item $g_1\sim g_2$ if and only if some positive powers of $g_1$ and
  $g_2$ are conjugate, and 
\item there exist rank 1 elements $g_1,g_2\in \Gamma$ such that
  $g_1\not\sim g_2$.
\end{enumerate}
\end{prop}

\begin{proof}
Analogous to the proof of Proposition 6 in \cite{bf}.
\end{proof}

\begin{thm}
Suppose that the action of $\Gamma$ on $X$ satisfies $WPD$. Then
$\widetilde {QH}(\Gamma)$ is infinite-dimensional.
\end{thm}

\begin{proof}
This follows from Proposition \ref{p6} and Theorem \ref{t1}.
\end{proof}

The Main Theorem, stated in the introduction, follows
(see Proposition \ref{DDFT}).

\section{Reducible nonpositively curved manifolds}\label{reducible}

In this section we will prove Theorem \ref{dichotomy}. We will need
three lemmas.

\begin{lemma}\label{A}
If $G$ is a group and $H\subset G$ has finite index, then
the restriction map $\widetilde{QH}(G)\to\widetilde{QH}(H)$ is
injective. In particular, if $\widetilde{QH}(H)=0$ then
$\widetilde{QH}(G)=0$ and if $\widetilde{QH}(G)$ is
infinite-dimensional then $\widetilde{QH}(H)$ is infinite-dimensional.
\end{lemma}

\begin{proof}
By passing to the intersection of conjugates of $H$, we may assume
that $H$ is a normal subgroup of $G$. We need to argue that if
$\phi:G\to\R$ is a homogeneous quasi-homomorphism such that $\phi|H$ is a
homomorphism, then $\phi$ is a homomorphism.

We first show that $\phi(gh)=\phi(g)+\phi(h)$ when $g\in G, h\in H$. Let $N$ be
the index of $H$ in $G$. Then we have:
$$
\phi(gh)=\frac 1N\phi((gh)^N)=\frac 1N\phi(ghg^{-1}\cdot g^2hg^{-2}\cdots
g^Nhg^{-N}\cdot g^N)
$$
with all expressions between the dots representing elements of $H$. Thus
$$
\begin{array}{l}
\phi(gh)=\frac
1N(\phi(ghg^{-1})+\phi(g^2hg^{-2})+\cdots+\phi(g^Nhg^{-N})+\phi(g^N))=\\
\frac 1N(N\phi(h)+N\phi(g))=\phi(h)+\phi(g)
\end{array}
$$
since $\phi$ is constant on conjugacy classes.

Now denote $H_0:=Ker(\phi|:H\to\R)$ and note that $H_0$ is normal in
$G$. Also note that by the above calculation $\phi$ is constant on the
cosets of $H_0$, so it induces a function $\overline
\phi:G/H_0\to\R$. This function is a homogeneous quasi-homomorphism on
the virtually abelian group $G/H_0$ (the subgroup $H/H_0$ is abelian
and has finite index). Therefore $\overline \phi$ is a homomorphism and
it follows that $\phi$ is a homomorphism.
\end{proof}

We now want to examine the other extreme, that is, whether having
infinite-dimensional $\widetilde{QH}$ is inherited by finite index
overgroups. In view of Lemma \ref{A} there is no harm in assuming
that the subgroup is normal, so we have an exact sequence
$$1\to H\to G\to\Sigma\to 1$$ with $\Sigma$ a finite group. There is a
natural action of $\Sigma$ on $\widetilde{QH}(H)$ (viewed as the space
of homogeneous quasi-homomorphisms on $H$ modulo $Hom(H,\R)$) given by
$$\sigma\cdot \phi(h)=\phi(\tilde\sigma^{-1} h \tilde\sigma)$$ for any
$\tilde\sigma\in G$ that maps to $\sigma$.

\begin{lemma}\label{C}
The injective map
$$\widetilde{QH}(G)\to\widetilde{QH}(H)$$ has image equal to 
$$\widetilde{QH}(H)^\Sigma$$ (the fixed subspace under the action of
$\Sigma$). 
\end{lemma}

\begin{proof}
If $\phi$ is a homogeneous quasi-homomorphism on $G$ then
$\phi(\tilde\sigma^{-1}h\tilde\sigma)=\phi(h)$, so $\phi|H$ is fixed
by $\Sigma$.

Conversely, suppose $\phi$ is a homogeneous quasi-homomorphism on $H$
which is fixed by $\Sigma$. We proceed as in the proof of Lemma
\ref{A}. Define $\tilde\phi:G\to\R$ by
$$\tilde\phi(g)=\frac 1N \phi(g^N)$$ where $N=|\Sigma|$. It is clear
that $\tilde\phi$ extends $\phi$ and we need to argue that
$\tilde\phi$ is a quasi-homomorphism. For $h\in H$ we have
$$\begin{array}{l}
\tilde\phi(gh)=\frac 1N \phi((gh)^N)=\frac 1N \phi(ghg^{-1}\cdot
g^2hg^{-2}\cdot g^Nhg^{-N}\cdot g^N)\sim\\ \frac 1N
(\phi(ghg^{-1})+\phi(g^2hg^{-2})+
\cdots+\phi(g^Nhg^{-N})+\phi(g^N))=\\
\frac
1N(N\phi(h)+N\tilde\phi(g))=\tilde\phi(h)+\tilde\phi(g)
\end{array}$$ where $\sim$
denotes bounded error independent of $g,h$. 
Similarly, $\tilde\phi(hg)\sim \tilde\phi(h)+\tilde\phi(g)$. 
Choose a set-theoretic section
$\Sigma\to G$, $\sigma\mapsto\tilde\sigma$, of $G\to\Sigma$.
Now if $g_1,g_2\in G$ then we can write $g_1=h_1\tilde\sigma_1$
and $g_2=\tilde\sigma_2 h_2$ for some $h_i\in H$. Then
$$\tilde \phi(g_1g_2)=\tilde\phi(h_1\tilde\sigma_1\cdot\tilde\sigma_2
h_2)\sim
\tilde\phi(h_1)+\tilde\phi(h_2)+\tilde\phi(\tilde\sigma_1\tilde\sigma_2)$$
while
$$\tilde\phi(g_i)\sim\tilde\phi(h_i)+\tilde\phi(\tilde\sigma_i)$$
Lemma now follows since there are only finitely many terms of the form
$\tilde\phi(\tilde\sigma_1\tilde\sigma_2)$ and $\tilde\phi(\tilde\sigma_i)$.
\end{proof}

It is not {\it a priori} clear why $\widetilde{QH}(H)^\Sigma$ couldn't
be finite-dimensional, or even 0, when $\widetilde{QH}(H)$ is
infinite-dimensional. In the following special case the fixed subspace
is infinite-dimensional.

\begin{lemma}\label{B}
Let $$1\to\Gamma_1\times\cdots\times\Gamma_k\to\Gamma\to\Sigma\to 1$$ be an
exact sequence with $\Sigma$ a finite group. Also assume
\begin{enumerate}[(1)]
\item Conjugation $Ad(\gamma)$ by any $\gamma\in\Gamma$ preserves the
  product structure on $\Gamma_1\times\cdots\times\Gamma_k$, i.e. for
  every $\gamma$ there is a permutation $\pi=\pi_\gamma\in S_k$ such
  that
$$Ad(\gamma)(\Gamma_i)=\Gamma_{\pi(i)}$$
\item If a factor is preserved by conjugation then the restriction is
  conjugation by an element of the factor, i.e. if $\pi_\gamma(i)=i$
  then there exists $\gamma_i\in\Gamma_i$ such that
  $$Ad(\gamma)(x)=\gamma_i x\gamma_i^{-1}$$
for every $x\in\Gamma_i$.
\end{enumerate}
If $\widetilde{QH}(\Gamma_1)$ is infinite-dimensional then
$\widetilde{QH}(\Gamma)$ is infinite-dimensional.
\end{lemma}

\begin{proof}
  Note that $$\widetilde{QH}(\Gamma_1\times\cdots\times\Gamma_k)=
  \widetilde{QH}(\Gamma_1)\times\cdots\times\widetilde{QH}(\Gamma_k)$$
  is also infinite-dimensional. In view of Lemma \ref{C}, to construct
  elements in this space that are fixed by $\Sigma$ we take $\phi\in
  \widetilde{QH}(\Gamma_1)$ (viewed as a quasi-homomorphism on
  $\Gamma_1\times\cdots\times\Gamma_k$ which is trivial on the other
  factors) and sum over the orbit. Let $\tilde
  \phi=\sum_{\sigma\in\Sigma}\sigma\cdot \phi$. By (1) and (2), the
  restriction of $\tilde \phi$ to $\Gamma_1$ equals $\phi$ multiplied
  by the number of elements of $\Sigma$ that preserve the first factor
  (this number is at least 1, since the identity preserves all
  factors). This (kind of a transfer) map sends an infinite linearly
  independent collection in $\widetilde{QH}(\Gamma_1)$ to an infinite
  linearly independent collection in
  $\widetilde{QH}(\Gamma_1\times\cdots\times\Gamma_k)^{\Sigma}$, which
  equals $\widetilde{QH}(\Gamma)$ by Lemma \ref{C}.
\end{proof}

\begin{proof}[Proof of Theorem \ref{dichotomy}]
Let $M'\to M$ be a finite cover such that $M'$ is a Riemannian product
$$M'=M_1\times\cdots\times M_k$$
with each $M_i$ of positive dimension and with $k\geq 1$ maximal
possible. By taking a further finite cover if necessary, we may assume
that each $M_i$ is a Riemannian topologically trivial fiber bundle 
$$T_i\hookrightarrow M_i\to B_i$$ where $T_i$ is a torus, $B_i$ is
non-positively curved with finite volume and the
universal cover $\tilde B_i$ of $B_i$ is either 
\begin{itemize}
\item a point,
\item a higher rank locally
symmetric space and the deck group is an irreducible lattice, or
\item the deck group contains rank 1 isometries.
\end{itemize}
Of course, $T_i$ could be a point. For details, see Eberlein's book
\cite[Section 7]{Eb}. Note that in general it is not possible to arrange that
$M_i$ is a Riemannian product $B_i\times T_i$. However, the universal
cover $\tilde M_i$ of $M_i$ is the Riemannian product $\tilde
B_i\times \tilde T_i$. In the rank 1 case this is the de Rham
decomposition of $\tilde M_i$ and in the higher rank case $\tilde B_i$
may further decompose as the product of symmetric spaces.

There are
now two cases.

{\bf Case 1.} Each $\tilde B_i$ is either a point or a higher rank
symmetric space. In this case we claim that
$\widetilde{QH}(\Gamma)=0$. Indeed, the Burger-Mozes theorem implies
that $\widetilde{QH}(\pi_1(B_i))=0$ for each $i$ and therefore
$$\widetilde{QH}(\pi_1(M_i))=
\widetilde{QH}(\pi_1(B_i))\times\widetilde{QH}(\pi_1(T_i))=0$$
and likewise $$\widetilde{QH}(\pi_1(M'))=0$$
Lemma \ref{A} now implies that $\widetilde{QH}(\Gamma)=0$.

{\bf Case 2.} There exists some $i$, say $i=1$, such that some deck
transformation of $\tilde B_i$ is 
a rank 1 isometry. In this case we will show that
$\dim\widetilde{QH}(\Gamma)=\infty$. 

The group $\Gamma$ acts by isometries on $$\tilde M=\tilde
B_1\times\cdots\times \tilde B_k\times E$$
where we collected all $\tilde T_i$ into the Euclidean factor
$E$. This action preserves the product structure in the following
sense. By a {\it slice of type $t$} we mean any set of the form 
$$b_1\times\cdots\times \tilde B_t\times\cdots\times
b_k\times e$$
for $b_i\in \tilde B_i$ and $e\in E$. Similarly, when
$T\subset\{1,\cdots,k\}$, we define slices of type $T$ by fixing
points in coordinates outside $T$.
\begin{enumerate}[(i)]
\item For any $\gamma\in\Gamma$ there is a permutation
  $\sigma_\gamma\in S_k$ such that $\gamma$ maps slices of type $t$ to
  slices of type $\sigma_\gamma(t)$, for $t=1,2,\cdots,k$, and
\item If $\sigma_\gamma(t)=t$ there is an isometry $\gamma_t:\tilde
  B_t\to\tilde B_t$ such that for every slice $S$ of type $t$ the diagram
$$
\begin{array}{ccc}
S&\overset\gamma{\longrightarrow}&\gamma(S)\\
p_t\downarrow&&\downarrow p_t\\
\tilde B_t&\overset{\gamma_t}\longrightarrow&\tilde B_t
\end{array}
$$
commutes, where $p_t$ is the projection to $\tilde B_t$.
\item More generally, if $T=\{t_1,\cdots,t_r\}$ is
  $\sigma_\gamma$-invariant, there is an isometry $\gamma_T:\tilde B_T\to
  \tilde B_T$, where $\tilde B_T=\tilde B_{t_1}\times\cdots\times
  \tilde B_{t_r}$ such that for every slice $S$ of type $T$ the
  diagram
$$
\begin{array}{ccc}
S&\overset\gamma{\longrightarrow}&\gamma(S)\\
p_T\downarrow&&\downarrow p_T\\
\tilde B_T&\overset{\gamma_T}\longrightarrow&\tilde B_T
\end{array}
$$
commutes, where $p_T$ is the projection to $\tilde B_T$.
\end{enumerate}

By $\Sigma\subset S_k$ denote the group of permutations
$\sigma_\gamma$ from (i) as $\gamma$ ranges over all elements of
$\Gamma$. Note that if $j_1$ and $j_2$ are in the same $\Sigma$-orbit
then $\tilde B_{j_1}$ and $\tilde B_{j_2}$ are isometric. Consider the
orbit of $1$, say $1,2,\cdots,r$. For $t\leq r$ let $\tilde\Gamma_t$
be the group of isometries of $\tilde B_t$ consisting of the
isometries $\gamma_t$ as in (ii) as $\gamma$ ranges over all elements
of $\Gamma$ with $\sigma_\gamma(t)=t$. Finally, let $\tilde\Gamma$ be
the group of isometries of $\tilde B_1\times\cdots\times\tilde B_r$
generated by the projections $\gamma_{\{1,\cdots,r\}}$ of elements of
$\Gamma$ (as in (iii)) and by
$\tilde\Gamma_1\times\cdots\times\tilde\Gamma_r$.

Then we have an exact sequence
$$1\to\tilde\Gamma_1\times\cdots\times\tilde\Gamma_r
\to\tilde\Gamma\to\tilde\Sigma\to 1$$ where $\tilde\Sigma$ is the
restriction of $\Sigma$ to $\{1,2,\cdots,r\}$. Note that (1) and (2)
of Lemma \ref{B} hold: this is obvious for elements of
$\tilde\Gamma_1\times\cdots\times\tilde\Gamma_r$, and if $\gamma\in\Gamma$
with $\sigma_\gamma(\{1,\cdots,r\})=\{1,\cdots,r\}$ and
$\sigma_\gamma(1)=1$, then (1) follows from (i) and for (2)
$Ad(\gamma_{\{1,\cdots,r\}})(x)=\gamma_1x\gamma_1^{-1}$ for
$x\in\tilde\Gamma_1$, by (ii)
(since $\gamma_1 \in \tilde \Gamma_1$ by definition).

Note that $\tilde\Gamma_j$ is discrete since it contains the deck
group $\pi_1(B_j)$ as a subgroup of finite index. Therefore, by the
Main Theorem, $\dim\widetilde{QH}(\tilde\Gamma_j)=\infty$, and
consequently by Lemma \ref{B}
$$\dim\widetilde{QH}(\tilde\Gamma)=\infty$$
The projection $\Gamma'$ of $\Gamma$ to $\tilde\Gamma$ has finite index (it
contains $\pi_1(B_1)\times\cdots\times\pi_1(B_r)$) and so by Lemma
\ref{A} we deduce $\dim\widetilde{QH}(\Gamma')=\infty$. Since
$\Gamma'$ is a quotient of $\Gamma$ it now follows that
$\dim\widetilde{QH}(\Gamma)=\infty$ (composing $\Gamma\to\Gamma'$ with
a homogeneous quasi-homomorphism $\Gamma'\to\R$ gives a homogeneos
quasi-homomorphism, so we have an infinite-dimensional space of
homogeneous quasi-homomorphisms on $\Gamma$; modding out
$Hom(\Gamma,\R)$ kills a finite dimensional subspace).
\end{proof}

\section{Teichmuller space with Weil-Petersson 
metric}\label{wp}

Let $\WP=\WP(S)$ be the Teichm\" uller space of a compact 
orientable surface $S$ of negative Euler characteristic with
the Weil-Petersson metric.
$\WP$ is a Riemannian manifold of negative sectional curvature, 
which is not complete, and any two points are joined by a (unique) Riemannian
geodesic. For basic properties of $\WP$ see Wolpert's survey paper \cite{Wo}.
A pseudo-Anosov element, $g$, in the mapping class
group $MCG(S)$ has an invariant geodesic, $l$, in $\WP$ called the axis
(see \cite{Wo}).
There exists a constant $R>0$, which depends on $g$, 
with the following property. Let $N_R(l)$ be the 
closed $R$-neighborhood of $l$. Then $g$ acts cocompactly on $N_R(l)$.
In this paper, we say $N_R(l)$ is {\it regular}
if this condition holds. The following result was proved by Jason
Behrstock \cite{behrstock} using different methods.

\begin{prop}\label{wp.rank1}
A pseudo-Anosov element $g$ is a rank-1 isometry
on $\WP$. 
\end{prop}

\begin{proof}
Let $l$ be the axis of $g$.
Take $R>0$ such that $N_R=N_R(l)$ is regular.  For $p,q \in \WP$, let
$p',q' \in l$ be the nearest points to $p,q$, respectively.  We
conctruct two ruled triangles (cf. Chap 4 in \cite{Ca}) from these
four points.  One triangle $\Delta$ is obtained by coning off the
geodesic $[p,q]$ with the cone point $p'$ by geodesics, and the other
one $\Delta'$ is by coning off the geodesic $[p',q']$ with the cone
point $q$.

The triangles $\Delta,\Delta'$ have a Riemannian metric 
induced from $\WP$. The metric is non-degenerate (unless a triangle
degenerates to a geodesic segment) since 
$\WP$ has negative sectional curvature.
Let $K$ denote the curvature of the induced metric 
on the triangles.
The curvature $K$ at each point in $\Delta, \Delta'$ is at most 
the upper bound of 
the sectional curvatures at the same point in $\WP$ (cf. \cite{Ca}).
All sides of the triangles are geodesics with respect to 
the induced metric since they are geodesics in $\WP$.
By the Gauss-Bonnet theorem, 
$$-\int_{\Delta} K ds = 
\pi - (a+b+c), -\int_{\Delta'} K ds = \pi -(a'+b'+c'),$$
where $a,b,c$ are the inner angles of $\Delta$, and $a',b',c'$ are the 
angles of $\Delta'$.
Therefore,
$$-\int_{\Delta} K ds  -\int_{\Delta'} K ds \le 2 \pi.$$

For each $0<r \le R$, let $N_r$ be the $r$-neighborhood
of $l$ in $\WP$, and $\partial N_r$ its boundary.
Now assume $d([p,q],l)> R$. 
We have $[p,q]\cap N_R =\emptyset$.
Then,
$c_r=\partial N_r \cap \Delta, c_r'=\partial N_r \cap \Delta'$
are curves embedded in $\Delta,\Delta'$, respectively.
The union $c_r \cup c_r'$ joins two 
points $c_r \cap [p,p']$ and  $c_r' \cap [q,q']$, whose
distance in $\WP$ is at least $|p'-q'|-2r$.
Therefore, $|c_r|+|c_r'| \ge |p'-q'|-2r$.
It follows that, by integrating $|c_r|+|c_r'|$ from $0$ to $R$, 
$$Area(N_R \cap \Delta) + Area(N_R \cap \Delta') \ge R(|p'-q'|-2R).$$

The sectional curvature of $\WP$ is negative, therefore $K<0$ on
$\Delta,\Delta'$ since they are ruled triangles.  Moreover, there
exists a constant $k>0$ such that the sectional curvature in $N_R$ is
at most $-k$ since $l$ is the axis of a pseudo-Anosov element $g$ and
$N_R$ is invariant by $g$.  It follows that on $N_R \cap \Delta, N_R
\cap \Delta'$ we have $K \le -k$ as well.  Therefore,
$$k (Area(N_R \cap \Delta) + Area(N_R \cap \Delta')) \le
-\int_{\Delta} K ds -\int_{\Delta'} K ds \le 2 \pi.$$ We obtain, if
$d([p,q],l) \ge R$,
$$R(|p'-q'| -2R) \le \frac{2\pi}{k},$$
so that 
$|p'-q'| \le \frac{2\pi}{kR}+ 2R.$

Now let $A$ be a ball disjoint from $l$. We will argue that the
projection of $A$ to $l$ has bounded size. If the radius $r_A$ of $A$ is
$\leq R$ the projection has diameter $\leq 2R$ since the projection is
distance non-increasing. Otherwise, let $A'$
denote the ball with the same center and with radius $r_A-R$. It
suffices to argue that the projection of $A'$ has bounded size. But
the argument above applies to the center $p$ and any point $q\in A'$,
and the Proposition follows.
\end{proof}

\begin{prop}\label{wp.wpd}
The action of $MCG(S)$ on $\WP$ is WPD.
\end{prop}

In the proof we need the following lemma. Let $g$ be a
pseudoAnosov element with axis $l$ and let $N_r(l)$ be regular.

\begin{lemma}\label{narrow}
For any $C>0$, there exist $L>0$ such that 
if a geodesic segment $[s,t]$ is contained in $N_C(l)$ and 
$|s-t| \ge L$, then there is a point $u \in [s,t]$ such 
that $u \in N_R(\gamma)$.
\end{lemma}

\begin{proof}
Suppose not.  Let $s',t' \in l$ be the projections of $s$ and $t$ to
$l$.  Let $\Delta$ be the ruled triangle with the base segment $[s,t]$
and the cone point $s'$, and $\Delta'$ the ruled triangle with the
base $[s',t']$ and the cone point $t$.  There is a constant $k>0$ such
that the sectional curvature in $N_R(\gamma)$ is at most $-k$,
therefore, the sectional curvature of the ruled triangles is at most
$-k$ on $\Delta \cap N_R(\gamma)$ and $\Delta' \cap N_R(\gamma)$.  By
an argument similar to the one for Proposition \ref{wp.rank1}, if
$[s,t] \cap N_R(l) =\emptyset$, then $|s-t| < L$, where $L$ depends
only on $R,k,C$, since $\Delta \cup \Delta'$ contains a rectangle of
width $r$ with one side $[s',t']$. The area of the rectangle has an
upper bound, therefore there is an upper bound for $|s'-t'|$, and for
$|s-t|$.
\end{proof}

Note that by convexity of the distance function, if $[s,t]$ is
contained in $N_C(l)$ then all points of $[s,t]$ except possibly for
initial and terminal segments of length $<L$ are contained in
$N_R(l)$.

\begin{proof}[Proof of Proposition \ref{wp.wpd}]
Let $g$ be a pseudo-Anosov map and let $c>0$. Set $C:=c+d(x_0,l)$, let
$N_R(l)$ be regular and
let $L$ be the number from Lemma \ref{narrow}. Choose $M$ so that
$|x_0-g^M(x_0)|>2L$. Now suppose that $|x_0-\gamma(x_0)|\leq c$ and
$|g^M(x_0)-\gamma g^M(x_0)|\leq c$. By the above remark, the midpoints
$p$ of $[x_0,g^M(x_0)]$ and $\gamma(p)$ of $[\gamma(x_0),\gamma
g^M(x_0)]$ are in $N_R(l)$, and by convexity $|p-\gamma(p)|\leq
c$. Thus the compact set $N_R(l)\cap N_c(p)$ contains
$\gamma(p)$. Since the action of the mapping class group on
Teichm\"uller space is discrete, there are only finitely many such
$\gamma$. 
\end{proof}

By Proposition \ref{wp.rank1} and \ref{wp.wpd},
Theorem \ref{mapping class group}
follows from the Main Theorem since $\WP$ is a CAT(0) space.

\bibliography{./ref}
%\begin{thebibliography}{99}

%\end{thebibliography}

\parbox{6cm}{
\noindent Mladen Bestvina

\noindent Mathematics Department

\noindent University of Utah

\noindent Salt Lake City, UT 84112, USA

\noindent {\tt bestvina@math.utah.edu}
}
\hfill
\parbox{5cm}{
\noindent Koji Fujiwara

\noindent Graduate School of Information Science

\noindent Tohoku University

\noindent Sendai, 980-8579, Japan

\noindent {\tt fujiwara@math.is.tohoku.ac.jp}
}

\end{document}